
\documentstyle{amsppt}
\magnification = \magstep 1
\pageheight{7 in}
\refstyle{A}

\topmatter
\title
Semistable 3-fold flips
\endtitle
\author
Alessio Corti
\endauthor
\address Department of Mathematics, University of Chicago, 5734 S.
University Avenue, Chicago IL 60637 \endaddress
\email corti\@ math.uchicago.edu \endemail
\date December 3, 1994 \enddate
\endtopmatter
\document

{\bf 1 Introduction.} In this short note, I will give a proof of the existence
of semistable 3-fold flips, which does not use the classification of
log terminal (i.e., quotient) surface singularities.
This permits us to avoid a calculation, which is a sort of logical
bottleneck in all the existing approaches to semistable 3-fold flips
([Ka1 pg. 158--159], [Sh pg. 386--389], [Ka3 pg. 483--486]), however different
they may look.
The message is that Shokurov's main reduction step, as refined in
[FA, Ch. 18], can be used profitably in the semistable case also. My main
motivation was to develop an approach to semistable flips that would have some
fighting chances in dimension 4, and the present paper is a first (small) step
in that direction.

We always work over an algebraically closed field of characteristic zero.

\definition{1.1 Definition} Let $X$ be a normal projective variety,
$B \subset X$ a ${\Bbb Q}$-Weil divisor such that $K_X+B$ is log terminal (this
notion is recalled in 2.2). Let $R \subset \overline {NE}(X)$ be an extremal
ray with $(K+B)\cdot R <0$, $\varphi_R : X \to U$ the contraction of $R$.
$\varphi_R$ is said to be a {\it flipping contraction} if the
$\varphi_R$-exceptional locus has codimension $\geq 2$.

A flip of $\varphi_R$ is by definition a variety $X^+$, together with a
morphism $\varphi^+: X^+ \to U$, such that $K^+ + B^+$ is
${\Bbb Q}$-Cartier and $\varphi^+$-ample.
\enddefinition

It is easy to see that the flip is unique if it exists. It is not so easy,
but true, that $K^+ + B^+$ is log terminal (see 2.4).
The flip conjecture asserts that
flips exist and that there is no infinite sequence of them.
An important special class of flips is that of semistable flips. These are the
flips that appear in the minimal model program for a semistable family.
In short, one is given the additional structure of a
projective morphism $f: X \to \Delta$, where $\Delta = \operatorname{Spec}
{\Cal O}$ is the spectrum of a discrete valuation ring ${\Cal O}$, with central
and generic points $0, \eta \in \Delta$. We denote $X_0$, $X_\eta$ the fibers
over $0$, $\eta$.
All extremal rays, contractions, flips, are compatible with this structure. The
starting point is a semistable family in the sense of Mumford, but the minimal
model program will soon introduce singularities, which we call semistable
terminal singularities (see 2.5). These are all terminal, but an important
point is that not all terminal singularities occur.

\definition{1.2 Definition} Let $f: X \to \Delta$ be as above,
$\varphi_R: X \to U$ a semistable flipping contraction. We say that
$\varphi_R$ is a {\it special} semistable flipping contraction if there
exists a component $S \subset X_0$ such that $S \cdot R \not = 0$.
\enddefinition

In this paper we prove the following:

\proclaim{1.3 Theorem} Semistable 3-fold flips exist if special
semistable 3-fold flips exist.
\endproclaim

This is satisfactory because special semistable flips are easier to construct.
Since I am unable to improve upon existing treatments, I will limit myself
to quickly
sketching two constructions. Both these constructions are based on the
following:

\proclaim{1.4 Lemma} Let $\varphi_R : X \to U$ be a special semistable 3-fold
flipping contraction. Then, in a neighborhood of every positive dimensional
fiber of $\varphi_R$, there exists a (reduced) surface $B \in |-K_X|$ having
Du Val singularities only.
\endproclaim

\demo{Proof} By assumption, there is a component $S_0$ of $S=X_0$ with
$S_0 \cdot R \not=0$. But $S \sim 0$, so in fact there are components
$S_1$, $S_2$ of $S$ with $S_1 \cdot R <0$ and $S_2 \cdot R >0$. By
[FA, 19.11], there is $B\in |-K-S|=|-K|$ such that $K+S+B$ is log canonical
(in the language of Shokurov and [FA], such a $B$ is called 1-complement).
But $S$ is a Cartier divisor, so $K+B$ is canonical, and so is $K_B=
K+B|_B$. \qed
\enddemo

Let ${\Cal H} = |-K_X|$. A way to formulate 1.4 is to say that $K+{\Cal H}$ is
canonical in a neighborhood of the flipping locus. I now outline the
two known methods to construct the flip:

{\bf 1.5} Use [Ka1, 8.7]: take a double cover $X^\prime \to X$ branched
along two general members of ${\Cal H}$. Then $X^\prime$ has canonical
singularities and $K_{X^\prime} \sim 0$. The flip of $X$ is a ${\Bbb Z}/2$
quotient of the flop of $X^\prime$. In [Ka1], canonical flops are
reduced to terminal flops
with the crepant descent method, generalized in [FA, Ch 6].

{\bf 1.6} This method is a private communication from S. Mori,
related to work of his student T. Hayakawa, and it allows, more generally,
to construct canonical $K+ {\Cal H}$-flops (these may be flips, flops or
inverse flips). One uses the crepant descent method in the form given
in [FA, Ch 6]. The proof is by induction on $e(X)$, the number of
$K+{\Cal H}$-crepant valuations, the starting point is a crepant divisorial
contraction $h:X_1 \to X$, then $e(X_1)=e(X)-1$ and by induction one is
reduced to the case $e=0$, where the flop is a genuine $K$-flop on
a 3-fold with terminal singularities of index 1. It seems that one needs
some knowledge about terminal singularities in order to construct the blow up
$h$.

In both methods, the flip is reduced to a flop. Flops are very difficult,
their existence rests upon: 1) the very hard
implication terminal Gorenstein $\Rightarrow$ isolated $cDV$
(this is not an issue in the semistable case) and 2) simultaneous
resolution of surface Du Val singularities.

It would be very interesting to find a better way to use the 1-complement
$B$ of 1.4.

\definition{1.7 Remark} If $\varphi : X \to U$ is a semistable
3-fold flip, and more generally for every 3-fold flip, $K+{\Cal H}$
is canonical in a neighborhood of the flipping locus.
In fact, Mori's proof of the existence of 3-fold flips consists
in establishing the above statement. We are, of
course, with our method, unable to see this even for semistable flips.
\enddefinition

After some preliminaries, the main result is proven in \S 4. In \S 5 I give
a short proof of the main point in [Sh]. The
material in the beginning of \S 3 was expanded, beyond what is
strictly required for the proof of theorem 1.3, to fit the needs of \S 5.

\definition{1.8 Acknowledgments} The material of this note was conceived
in July 1994, during the 3rd Utah Summer School on moduli of surfaces
of general type.
I thank the organizers for providing a very nice environment. I am also
very grateful to J. Koll\'ar and S. Mori for a very useful discussion
that solved two major technical problems.
\enddefinition

\bigskip
{\bf 2 Various kinds of singularities.}
We begin with the following:

\definition{2.1 Definition} Let $X$ be a smooth variety, $S\subset X$ be a
reduced Cartier divisor. We say that $S$ is a smooth normal crossing
divisor if

2.1.1 X is a curve, or:

2.1.2 Every irreducible component $S_0 \subset  S$ is smooth and, for each
such component, $(S-S_0)|_{S_0}$ is a smooth normal crossing divisor.
\enddefinition

We will work with the following definition of log terminal singularities:

\definition{2.2 Definition} Let $X$ be a normal variety, $B \subset X$ a ${\Bbb
Q}$-Weil
divisor. Assume that $B=\sum b_iB_i$ is a formal linear combination of reduced
and irreducible codimension 1 subvarieties $B_i \subset X$ with rational
coefficients $0<b_i\leq 1$.

The divisor $K+B$ is log terminal if it is
${\Bbb Q}$-Cartier and there exists a projective
morphism $h: Z \to X$, from a nonsingular variety $Z$, satisfying the following
conditions:

2.2.1 The exceptional locus $E\subset Z$ of $h$ is a divisor, and
$E\cup \operatorname{Supp} h^{-1}_\ast B$ is a divisor with smooth
normal crossings in $Z$.

2.2.2 Let $E=\sum E_i$ with $E_i$ irreducible. Then:
$$K_Z+ h^{-1}_\ast B +E=h^\ast (K_X+B)+ \sum a_i E_i$$
where $a_i>0$ for all $i$.
\enddefinition

In particular, if $K+B$ is log terminal, every component
$B_i \subset \lfloor B \rfloor$ is
normal and every intersection $B_i \cap B_j$ of two distinct such components
is irreducible. By [Sz], the notion just given corresponds to divisorial log
terminal
of [FA, 2.13.3], and is equivalent to weakly Kawamata log terminal of
[FA, 2.13.4] (see below 2.4). At the moment 2.2 is the best candidate for
the ``correct'' definition of log terminal singularities.

\definition{2.3 Definition} Let $X$ be a normal variety,
$B \subset X$ a ${\Bbb Q}$-Weil
divisor as in 2.2. The divisor $K+B$ is log canonical if it is
${\Bbb Q}$-Cartier and for all
morphisms $h: Z \to X$, from a nonsingular variety $Z$, with exceptional
divisors $E_i$:
$$K_Z+ h^{-1}_\ast B +E=h^\ast (K_X+B)+ \sum a_i E_i$$
with all $a_i \geq 0$.
\enddefinition

It is easy to see that log terminal singularities are preserved under
divisorial contractions. They are also preserved by flips, but this is
harder to see, due to an intrinsic drawback of definition 2.2: we
require the existence of a resolution satisfying some properties,
rather than asking a similar property of all resolutions as in 2.3.
This difficulty is resolved via the following
result due to Szab\'o [Sz]:

\proclaim{2.4 Resolution lemma} Let $X$ be an irreducible variety over an
algebraically closed field of characteristic 0, $S \subset X$ a
subvariety of pure codimension 1, $U \subset X$ a smooth open
subvariety such that $S\cap U \subset U$ is a smooth normal
crossing divisor. Then, there exists a projective
morphism $h:Z \to X$, from a smooth $Z$, satisfying the following conditions:

2.4.1 $h: h^{-1} U \to U$ is an isomorphism,

2.4.2 $h^{-1}\bigl(S \cup (X \setminus U)\bigr)$ is a smooth normal crossing
divisor in $Z$.\qed
\endproclaim

I now introduce a class of singularities large enough to allow the minimal
model program of a projective semistable family. For a smaller class in
dimension 3, see \S 5.

In what follows we fix a
projective morphism $f: X \to \Delta$, where $\Delta = \operatorname{Spec}
{\Cal O}$ is the spectrum of a discrete valuation ring ${\Cal O}$,
with central and generic points $0, \eta \in \Delta$.
We denote $S=X_0$, $X_\eta$ the fibers over $0$, $\eta$.
All minimal model programs,
divisorial contractions, flips, will be tacitly required
to be compatible with this structure.

\definition{2.5 Definition} Let $f: X \to \Delta$ be as above. We say that $f$
(or $X$,
when there is no danger of confusion) has {\it semistable} terminal
singularities if the following conditions are satisfied:

2.5.1 $X$ itself has terminal singularities,

2.5.2 $S=X_0=f^\ast (0)$ is reduced and $K_X+S$ is log terminal.
\enddefinition

In particular, a projective semistable family in the sense of Mumford has
semi\-stable terminal singularities.
Note that  every component $S_i$ of $S$ is normal and $K_{S_i}+
(S-S_i)|_{S_i}$ is log terminal, as a further consequence every intersection
$S_i \cap S_j$ of two distinct such components is irreducible and normal, etc.

An important observation is that $S \sim 0$ is linearly equivalent to 0, so
$K_X\sim K_X+S$ and a divisorial contraction (resp. flip) for $K_X$ is
the same thing as a divisorial contraction (resp. flip) for $K_X+S$. As a
consequence, semistable terminal singularities are preserved by
(semistable) divisorial contractions and flips.

I will now recall two results form classification theory. The first classifies
log terminal 3-fold singularities with ``large'' boundary divisor, up to
analytic equivalence. The result is an easy consequence of inversion of
adjunction [FA, 17.6] and [Sz], and is proven in
[FA, 16.15]:

\proclaim{2.6 Lemma} Let $x \in B \subset X$ be a 3-fold germ
with $\emptyset \not =B$ reduced and $K_X+B$ log terminal. Then $B$ has at
most 3 irreducible components, and:

2.6.1 If $B$ has three irreducible components, $x \in B \subset X$
is analytically isomorphic to: $0 \in (xyz =0) \subset {\Bbb A}^3$

2.6.2 If $B=B_1+B_2$ has two irreducible components, one of the following
happens:

2.6.2.1 $B_1$ and $B_2$ are both ${\Bbb Q}$-Cartier and
$x \in B \subset X$ is analytically isomorphic to:
$0 \in (xy =0) \subset {1 \over m}(q_1,q_2,1)$
where $(q_1, q_2, m)=1$.

2.6.2.1 Neither $B_1$ nor $B_2$ is ${\Bbb Q}$-Cartier and
$x \in B \subset X$ is analytically isomorphic to:
$0 \in (t=0) \subset \bigl(xy+tg(z, t)=0\bigr)$,
everything taking place in affine toric 4-space
${1 \over m}(q_1,q_2,1,a)$ with $(q_i, a, m)=(q_1, q_2, m)=1$.\qed
\endproclaim

The next result is a classification of semistable terminal singularities,
up to analytic equivalence. For a simple proof, see [Ka3, 4.1]:

\proclaim{2.7 Lemma}
Let $X$ be a 3-fold and $f: X \to \Delta$ have semistable terminal
singularities.
Let $t\in {\Cal O}$ be the uniformizing parameter, $x \in S=(t=0)$ be a
point, $r$ the index of $K_X$ at $x$.
Then $x \in X$ is analytically equivalent to one of the following:

2.7.1 $(xyz=t)\subset {\Bbb A}^4$

2.7.2 $(xy=t)\subset A$ where $A={1 \over r}(a, -a, 1, 0)$
for some $(a, r)=1$.

2.7.3 Two cases:

2.7.3.1 $r >1$ and
$\bigl((xy =f(z^r, t)\bigr)\subset A$
with $A$ as in 2.7.2. Here $(f(Z, t)=0)$ is an isolated curve singularity in
the $Z, t$-plane and $f(0, t)\not = 0$. Also $f(Z, 0)\not = 0$, otherwise
we are in case 2.7.2.

2.7.3.2 $r=1$ and $x \in X$ is an isolated singularity of the form:
$\bigl(g(x,y, z)=tf(x, y, z, t)\bigr)\subset {\Bbb C}^4$ where $\bigl(
(g(x, y, z)=0 \bigr)$ is a surface Du Val singularity and $f$ is  arbitrary.
\qed
\endproclaim

I wish to emphasize that 2.7 is quite easy, unlike the classification
of all terminal 3-fold singularities.

\bigskip
{\bf 3 The main construction.}

\definition{3.1 Definition} Let $X$ be a normal variety, $Y\subset X$ a closed
subvariety. We say that $X$ is ${\Bbb Q}$-factorial (resp. analytically,
resp. formally ${\Bbb Q}$-factorial) along $Y$ if every Weil divisor on the
Zariski germ
(resp analytic germ, resp. formal completion) of $X$ along $Y$ is
${\Bbb Q}$-Cartier. $X$ is (analytically, formally) ${\Bbb Q}$-factorial
if it is so along every subvariety $Y \subset X$ (it is clearly enough to
check this
at all closed points $y \in X$).
\enddefinition

If $X$ is ${\Bbb Q}$-factorial along $Y$, $X$ is ${\Bbb Q}$-factorial at
every point $y$ of $Y$. Not so (obviously) if $X$ is analytically or formally
${\Bbb Q}$-factorial along $Y$. The following however is easy:

\proclaim{3.2 Lemma} Let $X$ be analytically (resp.
formally) ${\Bbb Q}$-factorial along $Y$. If $Y \subset X$ has
codimension $\geq 2$, $X$ is analytically (resp. formally)
${\Bbb Q}$-factorial along every point $y$ of $Y$. \qed
\endproclaim

\definition{3.3 Definition} Let $\Delta$ be as usual, $X$ a normal variety
and $f: X \to \Delta$ a morphism, $Y \subset X$ a subvariety.
$X$ is {\it stably} (analytically, formally) ${\Bbb Q}$-factorial along $Y$
if, for every base change $\Delta^\prime \to \Delta$,
$X^\prime$ is (analytically, formally) ${\Bbb Q}$-factorial along $Y^\prime$,
where $X^\prime$ is the normalized pull-back, and $Y^\prime \subset X^\prime$
the inverse image.
\enddefinition

We now study stably analytically ${\Bbb Q}$-factorial semistable terminal
3-fold singularities. These are Kawamata's moderate singularities [Ka2]:

\definition{3.4 Definition} Let $X$ be a 3-fold, $f: X \to  \Delta$ a not
necessarily projective morphism. Let $t \in {\Cal O}$ be a parameter. $f$ has
{\it moderate singularities} if the analytic germ at every point $x \in X$ is
isomorphic to one of the following germs:

3.4.1 $(xyz=t)\subset {\Bbb C}^4$

3.4.2 $(xy=t)\subset A$ where $A={1 \over r}(a, r-a, 1, 0)$
for some $(a, r)=1$.

3.4.3 $(xy =z^r+ t^n)\subset A$
for some $n$, with $A$ as in 3.4.2.
\enddefinition

The following is proven in [Ka2] (and below):

\proclaim{3.5 Lemma} Let $X$ be a 3-fold, and let $f: X \to \Delta$ be a
projective morphism with semistable terminal singularities. There is
then a base change $\Delta^\prime \to \Delta$ and a small, not necessarily
projective morphism, $h:X^{\prime \prime}\to X^\prime$, where $X^\prime$ is the
pull-back, such that $X^{\prime \prime}$ has moderate
singularities ($h$ can be taken to be projective {\it locally analytically}
over $X^\prime$).
\endproclaim

\demo{Proof} I will give a quick sketch of the proof. Start with a singularity
of the form:
$$xy =f(z^r, t)$$
in affine toric 4-space ${1 \over r}(a, r-a, 1, 0)$. As in 2.7.3.1,
$f(Z, 0) \not=0$, and by the Weierstrass
preparation theorem there exists a
${\Bbb Z}/r$-equivariant analytic change of coordinates such that, in the
new coordinates, $f(Z, t)$ is a Weierstrass polynomial in $Z$:
$$f(Z, t)= \sum_{i=0}^k Z^i f_i(t)$$
where $f_i(t)$ is a convergent power series with $f_i(0)=0$. After a base
change $t=u^d$:
$$f(Z, u^d)= \prod_{i=1}^k Z-\varepsilon_i(u)$$
Using [Ko, 2.2], it is easy to construct a small projective partial resolution
with $k$ singular points:
$$xy=z^r-u^{n_i}$$
in $A$, where $n_i= \operatorname{ord} \varepsilon_i$. Taking $d$ large enough
one can do this on all of $X$ simultaneously, but we may loose projectivity.
One should also note that base changing introduces some quotient singularities
along the curves which are intersection of the components of the central
fiber. These are however easily resolved. \qed
\enddemo

The following is an immediate corollary of the proof just given:

\proclaim{3.6 Corollary} A semistable terminal 3-fold singularity is
stably analytically ${\Bbb Q}$-factorial if and only if it is moderate.\qed
\endproclaim

We now begin our proof of existence of semistable 3-fold flips.

\definition{3.7 Definition} Let $\varphi : X \to U$ be a semistable 3-fold
flip,
$C \subset X$ be the $\varphi$-exceptional set. The flip is said to be
moderate if $X$ has moderate singularities at every point of $C$.
\enddefinition

\proclaim{3.8 Lemma} Semistable 3-fold flips exist if moderate semistable
3-fold flips exist.
\endproclaim

The proof of 3.8 is based upon the following:

\proclaim{3.9 Lemma} Let $\varphi : X \to U$ be a moderate semistable flip.
Let $\varphi^+ : X^+ \to U$ be the flip of $\varphi$, $C^+$
the $\varphi^+$-exceptional set. Then $X^+$ has moderate
singularities at every point of $C^+$.
\endproclaim

\demo{Proof} By 3.6 $X$ has analytically ${\Bbb Q}$-factorial singularities.
It is well known then that $X^+$ has analytically ${\Bbb Q}$-factorial
singularities along $C^+$. Since $C^+ \subset X^+$ has
codimension $\geq 2$, by 3.2 $X^+$ has analytically ${\Bbb Q}$-factorial
singularities. This is true after base change because the base change of
the flip is the flip of the base change. So $X^+$ has stably analytically
${\Bbb Q}$-factorial singularities. Then $X^+$ has moderate singularities by
3.6. \qed
\enddemo

\demo{Proof of 3.8} This is standard using 3.9, let me give a quick outline.
Let $X \to U$ be a flip, $\Delta^\prime \to \Delta$ a base change
as in 3.5, $X^\prime \to U^\prime$ the base change $X^{\prime \prime} \to
X^\prime$ as in 3.5, with $X^{\prime \prime}$ moderate.
Note that $X^\prime$, $U^\prime$ are acted upon by the cyclic
group $G$ of the covering $\Delta^\prime \to \Delta$. As a first step
we run a minimal model program for $X^{\prime \prime}$ over $U^\prime$.
By 3.9, this consists of a finite number of moderate flips
$X^{\prime \prime} \dasharrow
X^{\prime \prime \prime}$. Let $X^{\prime +}$ be the relative canonical model
of $X^{\prime \prime \prime} \to U^\prime$, whose existence is granted by
the base point free theorem. Then $X^{\prime +}\to U^\prime$ is the flip
of $X^\prime \to U^\prime$, and $X^{\prime +}/G\to U^\prime /G=U$ is the
flip of $X \to U$.\qed
\enddemo

\noindent {\bf 3.10} From now on we fix a nonspecial moderate semistable
3-fold flip $\varphi:
X \to U$, $C \subset X$ the flipping material. We now begin the basic
construction for the proof of 1.3.

Let $m$ be very large and $H_0 \in |-mK_X|$ a smooth member. Let $\overline
H \subset U$ be a Cartier divisor satisfying the following conditions:

3.10.1 $\overline H$ contains $\overline H_0 =\varphi H_0$,

3.10.2 $K_U+ \overline H$ is log terminal outside $\varphi C$

(the existence of $\overline H$ is a consequence of the standard Bertini
theorem on the quasi projective variety $U \setminus \varphi C$).

The following is the main result of this section:

\proclaim{3.11 Lemma} Possibly after a base change, there exists
a projective morphism $h:Z \to X$ satisfying the following  conditions:

3.11.1 $Z$ is smooth and the $h$-exceptional set $E$ is a divisor.

3.11.2 $h:Z \setminus p^{-1}C \to X \setminus C$ is an isomorphism. In
particular $E \subset Z_0$ and $h: Z_\eta \to X_\eta$ is an isomorphism.

3.11.3 $Z_0$ is reduced and $Z_0 \cup h^{-1}_\ast H$ is a smooth normal
crossing divisor.

The following moreover is true:

3.11.4 For every birational morphism $g: Y \to U$, $N^1(Y/U)$ is generated by
the $g$-exceptional divisors and the components of $g^{-1}_\ast \overline H$
\endproclaim

\demo{Proof} By 3.10.2 and the resolution lemma 2.4, there is $h: Z \to X$
satisfying
3.11.1--3, with the possible exception that $Z_0$ may be nonreduced.
3.11.4 is also satisfied, because, as we will check momentarily, the
conditions of the following lemma 3.12 are met.
By 3.10.1 $\overline H$ contains $\overline H_0=\varphi H_0$, and
$\overline H_0$ is a generator of $WD(U)/CD(U)$ (this notation is introduced
in 3.12 below) because $X$ is ${\Bbb Q}$-factorial and $H_0$ is a generator
of $N^1(X/U)$. The conditions of 3.12 are therefore satisfied, and 3.11.4
holds.

We will achieve all properties after base change and semistable
reduction. Let $t$ be a parameter in $\Delta$, $\Delta^\prime \to \Delta$
a base change, $u$ a parameter in $\Delta^\prime$ and
$t=u^d$. Denote $X^\prime$, $Z^\prime$ the base change, $h^\prime$, $H^\prime
\subset X^\prime$ etc. the corresponding objects after base change.
By the semistable reduction theorem, if $d$ is divisible enough, there
is a projective resolution $h^{\prime \prime}:Z^{\prime \prime} \to Z^\prime$
such that $Z^{\prime \prime}_0$ is reduced and smooth normal crossing.
We will check that $Z^{\prime \prime} \to X^\prime$ satisfies all the required
properties. 3.11.4 is still true after base change, since $X$ is stably
${\Bbb Q}$-factorial, in fact this is the reason why we introduced
stable ${\Bbb Q}$-factorializations in 3.5 to begin with.
So we only need to show that $Z^{\prime \prime}_0\cup h^{\prime \prime -1}_\ast
H^\prime$ is a smooth normal crossing divisor, and in fact it is enough
to check this locally at every point. That this is the case is more or less
obvious, but we will try to explain it carefully.
To this end, we need to recall part
of Mumford's construction of the semistable reduction.

Locally analytically $Z={\Bbb A}^3$, and  $$Z_0 \cup h^{-1}_\ast H=
(\prod_1^k z_i^{n_i} \prod_{k+1}^l z_i=0)$$ where $z_i$ are
coordinates on ${\Bbb A}^3$,
$Z_0=(t=0)$, $t=\prod_1^k z_i^{n_i}$, and $h^{-1}_\ast H=(\prod_{k+1}^l z_i=0
)$. After the base change $t=u^d$, the fiber product $Z^\prime$
is described as:
$$Z^\prime=A \times {\Bbb A}^{3-k}$$ where $A=(u^d=\prod_1^k z_i^{n_i})
\subset {\Bbb A}^{k+1}$, and, inside $Z^\prime$, $Z^\prime_0=(u=0)$ and
$h^{\prime -1}_\ast H^\prime = (\prod_{k+1}^l z_i=0)$ in coordinates
$z_{k+1}, ... z_3$ for ${\Bbb A}^{3-k}$. The construction of the semistable
reduction begins with taking the normalization of $A$:
$$\nu: A^\nu =\coprod_{j=1}^e A_j \to A$$
Here $(e)=(d, n_i)$ and each $A_j$ is isomorphic to the simplicial
affine toric variety ${1\over d}(n_1, ...n_k)$. Now $Z^{\prime \prime}$ is
obtained by gluing together pieces of the form:
$$B_j \times {\Bbb A}^{3-k} \to A_j \times {\Bbb A}^{3-k}$$
where $B_j \to A_j$ is a toric resolution of $A_j$ with the
property that $(u=0)
\subset B_j$ is a smooth normal crossing divisor (the proof of the
semistable reduction theorem consists in proving that
such $B_j \to A_j$ exist, and showing that a choice exists
so that the gluing is possible). Now $h^{\prime \prime}_\ast
H^\prime$ is described, in $Z^{\prime \prime}$, by the equation
$\prod_{k+1}^l z_i=0$, which proves the statement. \qed
\enddemo

\proclaim{3.12 Lemma} Let $U$ be a normal variety, $D_1, ..., D_k \subset U$
Weil
divisors on $U$ generating:
$${WD(U)\over CD(U)} =  {   \{{\Bbb Q} \text{-Weil divisors
on}\; U\} \over \{{\Bbb Q}\text{-Cartier divisors on}\; U\}}$$
Let $f: Y\to U$ be a birational
morphism. Then $N^1(Y/U)$ is generated by the $f^{-1}_\ast D_i$ and the
irreducible $f$-exceptional divisors.
\endproclaim

\demo{Proof} Let $A \in N^1(Y/U)$. A linear combination
$\sum \lambda_i D_i +f_\ast A$ is ${\Bbb Q}$-Cartier on $U$. Then
$\sum \lambda_i f^{-1}_\ast D_i +A -f^\ast (\sum \lambda_i D_i +A)$ is
$f$-exceptional. But in $N^1(Y/U)$, $f^\ast (\sum \lambda_i D_i +A)=0$. \qed
\enddemo

\bigskip
{\bf 4 Proof of 1.3.} We use Shokurov addition and subtraction
method with the refinements in
[FA 18.12]. We use the notation of \S 3, with one change: we now denote
$H$ the divisor $h^{-1}_\ast H$ on Z.
Let $\varphi \circ h = p:Z \to U$. Before coming to the
details, I will explain the broad outline of
the proof, in three steps as follows:

{\bf A} We run a MMP for the divisor $K_Z+H\sim K_Z+H+Z_0$, over
$U$. We need to show that flips exist and that each step of the
program preserves condition 2.5.1, namely all varieties involved in the
program have terminal singularities (all steps clearly preserve 2.5.2).
The program terminates at the
variety $p^\prime: Z^\prime \to U$ and $K_{Z^\prime}+H^\prime$ is
$p^\prime$-nef.

{\bf B} Using the nef threshold method, we progressively subtract pieces of
$H^\prime$ until we reach a variety $p^{\prime \prime}: Z^{\prime \prime}
\to U$ where $K_{Z^{\prime \prime }}$ is $p^{\prime \prime}$-nef.
As in A, we need to show that flips exist and that each step of the
program preserves condition 2.5.1.

{\bf C} The relative canonical model $\varphi ^+ : X^+ \to U$ of $p^{\prime
\prime} : X^{\prime \prime} \to U$, whose existence is granted by the
base point free theorem, is the flip of $\varphi$. This is obvious.

We will use the following notation:

{\bf 4.1} The central fiber $$Z_0=p^\ast U_0= \sum_{i=0}^I S_i$$
The $S_i$ for $i>0$ are irreducible components and
$S_0=p^{-1}_\ast U_0$. The $S_i$ for
$i \geq 1$ are precisely the $p$-exceptional components.
Note that the image of every
$S_i$ lies in $\overline H$, which follows from 3.11.2 and the fact that
$\overline H \supset \varphi C$

We now carry out steps A and B in detail:

{\bf A} The MMP for $K_Z+H\sim K_Z+H+Z_0$ constructs varieties
$Z=Z^1 \dasharrow Z^2 \dasharrow \cdots Z^\alpha$, with projections
$p^\alpha : Z^\alpha \to U$, divisors $Z^\alpha_0=
\sum_{i=0}^{I(\alpha)}S_i^\alpha$, etc..
Let us inductively assume that
$Z^\alpha$ has already
been constructed, and that property 2.5.1 holds for $Z^\alpha$. Let
$\psi^\alpha : Z^\alpha \to V^\alpha$ be an extremal contraction,
corresponding to a ray
$R^\alpha$ such that $(K_{Z^\alpha}+H^\alpha)\cdot R^\alpha <0$. I will in
the sequel drop the superscript $\alpha$ from the notation.

{\bf A.1} If $\psi$ is a divisorial contraction, we need to show that
property 2.5.1 holds for $V$. Since $H$ contains no proper divisor, $H
\cdot R \geq 0$.
The ray is a ray for $K_Z$ so 2.5.1 is clearly preserved.

{\bf A.2} If $\psi$ flips, we need to show that the flip $\psi^+ : Z^+ \to V$
exists and that property 2.5.1 holds for $Z^+$.
Let $C$ be a connected component of the flipping set. We distinguish two
cases, according to wether there exists an irreducible component $M$ of
$H$ with $M \cdot C <0$ or not:

{\bf A.2.1} There is an irreducible component  $M$ of
$H$ with $M \cdot C <0$. As $C$ is contained in $S_i$ for some $i$, by
2.6.2.1 we must have that $M$, $S_i$ are irreducible in a neighborhood of $C$,
and $C=M\cdot S_i$. Because $C \subset M$ is contractible, $S_i \cdot C <0$
and there exists $S_j$ with $S_j \cdot C >0$. By 2.6 again, and since by
assumption $(K+H+Z_0)\cdot C <0$, $C$ intercepts no component of $H$
other than $M$. Then $K+H$ satisfies condition $(\ast)$ of [FA, 5.1] in a
neighborhood of $C$. Incidentally, note that $\psi : M \to \psi(M)$ is a
2-dimensional semistable extremal contraction (for $K_M \sim K_M +
M \cdot Z_0$).
The flip exists by [FA, 5.4.1], and condition $(\ast)$ holds for $Z^+$ by
[FA, 5.2]. In particular 2.5.1 holds for $Z^+$.

{\bf A.2.2} $M\cdot C \geq 0$ for all irreducible components $M$ of
$H$. In particular the flip is also a flip for $K_Z=K_Z+Z_0$ so $Z^+$
satisfies 2.5.1 if it exists. Let us prove that $Z^+$ exists. By
3.11.4 there is a component $M$ of
$H+\sum_{i>0}S_i$ such that $M \cdot C \not = 0$. Now:
$$ p^\ast \overline H = H + \sum_{i>0} b_i S_i$$
and all $b_i >0$, since we made sure that the $p$-image of
every $p$-exceptional divisor lies in $\overline H$. Since $p^\ast
\overline H \cdot C =\overline H \cdot p_\ast C=0$, there is a component
$L$ of $H+\sum_{i>0}S_i$ such that $L \cdot C <0$. By assumption, $L$ is
one of the $S_i$. The flip is a special flip, and it exists by assumption.

{\bf B} Let $p^\prime: Z^\prime \to U$ with $K_{Z^\prime}+H^\prime$
$p^\prime$-nef be the model constructed in A. The nef threshold method is a
guided version of the minimal model program.
Without attempting a general formulation, I will describe what this is in
the present circumstances.
We create varieties $Z^{\prime}=Z^{\prime 1} \dasharrow Z^{\prime 2} \cdots
\dasharrow Z^{\prime \alpha}$ and rational numbers $\varepsilon^\alpha$
such that $K_{Z^{\prime \alpha}}+\varepsilon^\alpha H^{\prime \alpha}$ is nef.
Assume $Z^{\prime \alpha}$ has been constructed, I will now give the recipe
for $Z^{\prime \alpha +1}$.
$\varepsilon^{\alpha +1}$ is the smallest $\varepsilon \leq
\varepsilon^\alpha$ such that $K_{Z^{\prime \alpha}}+\varepsilon
H^{\prime \alpha}$ is still nef, and $Z^{\prime \alpha} \dasharrow
Z^{\prime \alpha + 1}$ the modification associated to one of the (finitely
many) extremal rays $R^\alpha$ such that $$\bigl(K_{Z^{\prime
\alpha}}+(\varepsilon^{\alpha +1}-\eta)
H^{\prime \alpha} \bigr) \cdot R^\alpha <0$$
for all $\eta$ very small.

Assuming that the contraction $\psi_{R^\alpha} :
Z^\alpha \to V^\alpha$ is flipping, I need
to prove that the flip exists.
I will in
the sequel drop the superscript $\alpha$ from the notation.
By construction $H^\prime \cdot R >0$. But $H^\prime \cdot R = -
\sum b_i S^\prime_i \cdot R$,
so $S_i^\prime \cdot R < 0$ for some $i$.
As $Z^\prime_0 \cdot R =0$, $S_j^\prime \cdot R
>0$ for some $j$. The flip is special, and it exists by assumption.
This concludes the proof of 1.3. \qed

\definition{4.2 Remark} As a final remark, I wish to say that, using
the $G$-minimal model program and a $G$-invariant version of Mumford
semistable reduction, it would have been possible to prove 1.3 avoiding
the material on stable ${\Bbb Q}$-factorialization and the classification
2.7 (but not 2.6) altogether. Such an argument would have probably been
more complicated than the one given in the text.
\enddefinition

\bigskip
{\bf 5 Strictly semistable 3-fold singularities.} In this last section, I will
prove some results in [Sh]. I would be interested in knowing
if the material here has some sort of higher dimensional generalization.

In this section all varieties, germs, etc. are tacitly assumed to come
equipped with a morphism, usually a projective one, to $\Delta$.
All minimal models, divisorial cotractions, etc., are required to be compatible
with this structure.

\definition{5.1 Definition} A germ $x \in X$ (resp. a variety $X$) admits a
semistable resolution if there is a resolution $f: Z \to X$ such that
$Z_0 \subset Z$ is a reduced smooth normal crossing divisor (in other
words, $Z$ is semistable in the sense of Mumford).
\enddefinition

The following is well known:

\proclaim{5.2 Lemma} Let $x \in X$ be a moderate 3-fold singularity.
Then $x \in X$ admits a projective semistable resolution.
\endproclaim

\demo{Proof} Let $x \in X$ be described by the equation:

5.2.1 $(xy=t)\subset A$, or

5.2.2 $(xy =z^r+ t^n)\subset A$,

where $A={1 \over r}(a, r-a, 1, 0)$ for some $(a, r)=1$. In both cases
let $f: B \to A$ be the weighted blow up with weights
${1 \over r}(a, r-a, ,1, r)$, $X^\prime \subset B$ the proper preimage,
and let us also denote $f:X^\prime \to X$ the restriction to $X^\prime$.
The following are easily checked via a small calculation in explicit
coordinates:

In case 5.2.1, $X^\prime$ has two singular points, given by
$(xy=t)\subset {1 \over a}([r], [a-r], 1, 0)$ and
$(xy=t)\subset {1 \over r-a}([r], [a-r], 1, 0)$

In case 5.2.2, $X^\prime$ has three singular points, given by
$(xy=t)\subset {1 \over a}([r], [a-r], 1, 0)$,
$(xy=t)\subset {1 \over r-a}([r], [a-r], 1, 0)$,
and $(xy=z^r+t^{n-1})\subset A$.

It is then immediate that a repeated application of
$f: X^\prime \to X$ gives the desired resolution. \qed
\enddemo

\definition{5.3 Definition} A semistable terminal 3-fold singularity is
{\it strictly} semistable if its analytic ${\Bbb Q}$-factorialization is
moderate.
\enddefinition

As an immediate consequence of 5.2, we have:

\proclaim{5.4 Corollary}

5.4.1 A strictly semistable terminal analytic germ admits a projective
semistable resolution.

5.4.2 If $X$ has strictly semistable terminal singularities, there is
a semistable resolution $Z \to X$, which is projective locally analytically
in $X$.
\endproclaim

By 3.9 strictly semistable terminal singularities are preserved by flips.
Our goals in this section are to show that they are also preserved
by divisorial
contractions, and to establish a converse to 5.4. The first result
is due to Shokurov [Sh].

\proclaim{5.5 Theorem} Assume that $X$ is projective over $\Delta$ and has
strictly semistable terminal singularities. Let $f: X \to Y$ be an extremal
divisorial contraction. Then $Y$ has strictly semistable terminal
singularities.
\endproclaim

\proclaim{5.6 Theorem} Let $x \in X$ be a semistable terminal analytic germ.
The following are quivalent:

5.6.1 $x \in X$ is strictly semistable,

5.6.2 $x \in X$ has a projective semistable resolution.
\endproclaim

5.6 is an easy consequence of 5.5:

\demo{Proof of 5.6 using 5.5} We need to show that 5.6.2 implies 5.6.1.
Let $Z \to X$ be a projective semistable resolution. Run a minimal model
program for $Z \to X$. By 3.2, this terminates at an analytic
${\Bbb Q}$-factorialization $X^\prime \to X$. By 5.5, $X^\prime$ has strictly
semistable terminal singularities, and so does $x\in X$. \qed
\enddemo

The proof of 5.5 uses the following:

\proclaim{5.7 Lemma} Let $x \in X$ be an analytically ${\Bbb Q}$-factorial
rational singularity. If $x \in X$ admits a semistable resolution, it
is stably ${\Bbb Q}$-factorial. In particular, if $x\in X$ is semistable
terminal, it is moderate.
\endproclaim

\demo{Proof of 5.5 using 5.7} Let $E \subset X$ be the exceptional divisor.
Abusing notation I will in the sequel denote $X$ the analytic germ of
$X$ along $E$. Let $y=f(E)$. The point here is that $Y$ is not necessarily
analytically ${\Bbb Q}$-factorial along $y$, and $X$ is not necessarily
analytically ${\Bbb Q}$-factorial along $E$. Let $X^\prime \to X$ be
a projective analytic ${\Bbb Q}$-factorialization along $E$. Run now
a minimal model program for $X^\prime \to Y$. After a finite number
of flips $X^\prime \dasharrow X^{\prime \prime}$ I meet a divisorial
contraction $X^{\prime \prime} \to Y^\prime$. Here $Y^\prime \to Y$ is an
analytic ${\Bbb
Q}$-factorialization. Note that, since flips preserve strictly semistable
terminal singularities, $X^{\prime \prime}$ has strictly semistable
singularities, so it admits a semistable resolution $Z \to X^{\prime \prime}$
Then the composition $Z \to Y^\prime$ is a semistable resolution, and
since $Y^\prime$ is analytically ${\Bbb Q}$-factorial, by 5.7 it is stably
analytically ${\Bbb Q}$-factorial, hence moderate. This proves 5.5. \qed
\enddemo

\demo{Proof of 5.7} Let $\Delta^\prime \to
\Delta$ a base change of degree $d$.
Let $X^\prime \to \Delta^\prime$ be the normalization of the fiber product,
and $\pi: X^\prime \to X$ the natural map. Fix a semistable resolution
$h: Z \to X$ and let $Z^\prime$ be the fiber product.
The following are the two crucial observations:

5.7.1 $Z^\prime$ is already normal, and in particular there is a natural
map $h^\prime : Z^\prime \to X^\prime$, given by the Stein factorization.
More to the point, $Z^\prime$ has toroidal simplicial, hence analytically
${\Bbb Q}$-factorial, singularities.

5.7.2 The cyclic group $G=\mu_d$ acts on $Z^\prime$, $X^\prime$ in such a way
that $Z^\prime/G =Z$, $X^\prime/G=X$ and $h^\prime : Z^\prime
\to X^\prime$ is $G$-equivariant. Most importantly, $G$ acts trivially on the
central fiber $Z^\prime_0$.

Let $D^\prime \subset X^\prime$ be a Weil divisor. Our aim is to show that
$D^\prime$ is ${\Bbb Q}$-Cartier. Certainly $\sum_{g\in G} gD^\prime
=\pi^\ast \pi_\ast D^\prime$ is ${\Bbb Q}$-Cartier, and so is its
pull back $h^{\prime \ast}\sum gD^\prime =\sum gh^{\prime -1}_\ast D^\prime
+ \sum da_iE_i$, where $E_i \subset Z^\prime_0$ are the $h^\prime$-exceptional
components. First of all I claim that:
$$h^{\prime -1}_\ast D^\prime + \sum a_iE_i \equiv 0$$
is numerically equivalent to zero relatively to $h^\prime$. Let indeed $C
\subset Z^\prime_0$ be
a curve such that $h^\prime C$ is a point. Then, using 5.7.1
to intersect with $h^{\prime -1}_\ast D^\prime$:
$$h^{\prime -1}_\ast D^\prime
\cdot C=gh^{\prime -1}_\ast D^\prime
\cdot gC =gh^{\prime -1}_\ast D^\prime
\cdot C$$ by 5.7.2. Then:
$$0=(\sum gh^{\prime -1}_\ast D^\prime
+ \sum da_iE_i)\cdot C=
d(h^{\prime -1}_\ast D^\prime + \sum a_iE_i)\cdot C$$
which shows the claim. $X$ has rational singularities, and so does $X^\prime$.
$R^1h^\prime_\ast {\Cal O}_{Z^\prime}=(0)$, hence a multiple $\nu (h^{\prime
-1}_\ast D^\prime + \sum a_iE_i)\sim {\Cal O}_{Z^\prime}$ is {\it linearly}
equivalent to zero.
Then $\nu D^\prime$ is Cartier. \qed
\enddemo

\enddocument